\documentclass[11pt]{amsart}

\theoremstyle{plain}
\newtheorem{theorem}                 {Theorem}      [section]
\newtheorem{proposition}  [theorem]  {Proposition}
\newtheorem{corollary}    [theorem]  {Corollary}

\theoremstyle{definition}
\newtheorem{example}      [theorem]  {Example}
\newtheorem{remark}       [theorem]  {Remark}

\numberwithin{equation}{section}

\def \1{\mbox{${\mathbf 1}$}}
\def \r{\mathbb R}
\def \s{\mbox{${\mathbb S}$}}

\def \t{\mbox{${\partial/\partial t}$}}
\def \ss{\mbox{${\partial/\partial s}$}}
\def \tt{\mbox{${\scriptstyle\partial/\partial t}$}}

\def \F{\mbox{$\phi$}}

\def \ts{\mbox{${\tiny \s}$}}

\def \rrm{\mbox{${\scriptstyle \frac{1}{\sqrt 2}}$}}

\DeclareMathOperator{\trace}{trace}
\DeclareMathOperator{\grad}{grad}

\DeclareMathOperator{\riem}{Riem}
\DeclareMathOperator{\ricci}{Ricci}

 \DeclareMathOperator{\Hom}{Hom}

\begin{document}
\title{Submanifolds with biharmonic Gauss map}
\author{A. Balmu\c s}
\address{Faculty of Mathematics, ``Al.I.~Cuza'' University of Iasi\\
\newline
Bd. Carol I Nr. 11 \\
700506 Iasi, ROMANIA} \email{adina.balmus@uaic.ro}
\author{S. Montaldo}
\address{Universit\`a degli Studi di Cagliari\\
Dipartimento di Matematica\\
\newline
Via Ospedale 72\\
09124 Cagliari, ITALIA} \email{montaldo@unica.it}
\author{C. Oniciuc}
\address{Faculty of Mathematics, ``Al.I.~Cuza'' University of Iasi\\
\newline
Bd. Carol I Nr. 11 \\
700506 Iasi, ROMANIA} \email{oniciucc@uaic.ro}

\subjclass[2000]{58E20}

\thanks{}

\keywords{Biharmonic maps, Gauss map, hypercones, isoparametric
hypersurfaces.}

\maketitle

\begin{abstract}
We generalize the Ruh-Vilms problem by characterizing the
sub\-ma\-ni\-folds in Euclidean spaces with proper biharmonic Gauss
map and we construct examples of such hypersurfaces.
\end{abstract}

\section{Introduction}
As it is classically known, most of the extrinsic geometry of an
oriented submanifold $M^m$ in the Euclidean space $\r^{m+n}$ can be
described by its Gauss map $\gamma:M\to G(m,n)$ which  assigns to
every point $p\in M$ the tangent space $T_pM$, thought of as a point
of the Grassmannian of oriented $m$-dimensional subspaces of
$\r^{m+n}$. A splendid example is the celebrated Ruh-Vilms Theorem
which asserts that the Gauss map $\gamma:M\to G(m,n)$ is a harmonic
map if and only if the mean curvature vector field of $M$ in
$\r^{m+n}$ is parallel. Here we say that a smooth map $\phi:(M,g)\to
(N,h)$ between Riemannian manifolds is {\it harmonic} if it is a
critical point of the {\it energy functional}
$E(\phi)=\frac{1}{2}\int_{M}\, |d\phi|^2\,v_g$, i.e. $\phi$ is a
solution of the  corresponding Euler-Lagrange equation which is
given by the vanishing of the tension field $\tau(\phi)=\trace\nabla
d\phi$.

A natural extension of harmonic maps is provided by {\it biharmonic
maps} (as suggested by J.~Eells and J.H.~Sampson in \cite{ES}) which
are the critical points of the {\it bienergy functional}
$E_2(\phi)=\frac{1}{2}\int_{M}\, |\tau(\phi)|^2\,v_g$. In \cite{Jia}
G.Y.~Jiang derived the first variation formula of the bienergy
showing that the Euler-Lagrange equation for $E_2$ is
\begin{eqnarray}\label{eq:bihar1}
\tau_2(\phi)&=&-J(\tau(\phi))=
-\Delta\tau(\phi)-\trace R^{N}(d\phi,\tau(\phi))d\phi\notag \\
&=&0,
\end{eqnarray}
where $J$ is (formally) the Jacobi operator of $\phi$, $\Delta$ is
the rough Laplacian defined on sections of $\phi^{-1}(TN)$ and
$R^{N}(X,Y)=\nabla_X\nabla_Y-\nabla_Y\nabla_X-\nabla_{[X,Y]}$ is the
curvature operator on $(N,h)$.

In this paper we propose to study the {\it biharmonic equation}
($\tau_2(\phi)=0$) for the Gauss map of submanifolds in the
Euclidean space, in the intent to generalize the Ruh-Vilms Theorem
to the case of biharmonicity. To pursue our intent  we first derive
the equation that characterizes the submanifolds in the Euclidean
space with biharmonic Gauss map (Theorem~\ref{th: bih_gauss}).
Although the condition that ensures the biharmonicity of the Gauss
map is rather technical, in the case of hypersurfaces it simplifies
and gives the following: the Gauss map of an orientable hypersurface
$M^m$ in $\r^{m+1}$ is proper biharmonic if and only if $\grad f\neq
0$ and
$$
\Delta\grad f +A^2(\grad f)-|A|^2\grad f=0,
$$
where $\Delta$ denotes the rough Laplacian on $C(TM)$ while $f$ and
$A$ denote the mean curvature function and the shape operator,
respectively (see also~\cite{AB}).

The last part of the paper is devoted to the construction of
examples of hypersurfaces with biharmonic Gauss map. We study the
biharmonicity of the Gauss map for hypercones generated by constant
mean curvature hypersurfaces in spheres (Theorem~\ref{th:
CMC_hypercone}) and, in particular, by isoparametric hypersurfaces,
obtaining explicit  examples. Non-existence results for hypercones
in $\r^3$ and $\r^4$ with proper biharmonic Gauss map are obtained
(Theorem~\ref{th: non-ex_cone_Gauss_R3} and Theorem~\ref{th:
non-ex_cone_Gauss_R4}).

\section{Preliminaries}

\subsection{Biharmonic maps between Riemannian manifolds}

We recall the following facts on biharmonic maps:
\begin{enumerate}
\item[(i)] the equation $\tau_2(\phi)=0$ is called the {\it biharmonic equation}
and a map $\phi$ is biharmonic if and only if its tension field is
in the kernel of the Jacobi operator;

\item[(ii)] a harmonic map is obviously a biharmonic
map. We call {\it proper biharmonic} the biharmonic non-harmonic
maps;

\item[(iii)] a harmonic map is an absolute minimum of the
bienergy;

\item[(iv)] if $M$ is compact and $\riem^{N}\leq
0$, i.e. the sectional curvature of $(N,h)$ is non-positive, then
$\phi:M\to N$ is biharmonic if and only if it is harmonic;

\item[(v)] if $\phi:M\to N$ is a Riemannian
immersion with $|\tau(\phi)|=$ constant and $\riem^{N}\leq 0$, then
$\phi$ is biharmonic if and only if it is harmonic (minimal).
\end{enumerate}
The first three remarks are immediate consequences of the definition
of the bienergy and of \eqref{eq:bihar1}. The non-existence results
(iv) and (v) are proved in \cite{J1} and in \cite{O2}, respectively.

\noindent On the other hand, in Euclidean spheres we do have
examples of proper biharmonic submanifolds, i.e. non-minimal
submanifolds for which the inclusion map is biharmonic. It was
conjectured in ~\cite{ABSMCO1} that the only proper biharmonic
hypersurfaces in $\s^{m+1}$ are the open parts of the hypersphere
$\s^m(\rrm)$ and of the generalized Clifford torus
$\s^{m_1}(\rrm)\times\s^{m_2}(\rrm)$, $m_1+m_2=m$ and $m_1\neq m_2$.

For a general account on biharmonic maps see ~\cite{SMCO}.

\subsection{The Gauss map}
Consider $M^m$ to be a $m$-dimensional oriented submanifold in
$\r^{m+n}$. The map which assigns to every point $p\in M$ the
oriented tangent space $T_pM$, thought of as a point of the
Grassmannian of oriented $m$-dimensional subspaces of $\r^{m+n}$,
\begin{align*}
\gamma:M&\longrightarrow G(m,n)\\
p&\longmapsto T_pM,
\end{align*}
is called the {\it Gauss map associated to $M$}.

As usually, the Riemannian structure on $G(m,n)$ is defined by
considering the Euclidean metric on $\r^{m+n}$ and by identifying
the tangent space to $G(m,n)$ at a point $P\in G(m,n)$ as follows:
$$
T_PG(m,n)=\Hom(P,P^\perp)=P^*\otimes P^\perp.
$$
Thus, if we fix a positive oriented orthonormal basis
$e_1,\ldots,e_m$ of $P$ and complete it to an orthonormal basis of
$\r^{m+n}$ with $e_{m+1},\ldots, e_{m+n}$, spanning $P^\perp$, then
a basis of $T_PG(m,n)$ will be given by
$$\{e^*_i\otimes
e_{m+a}\}_{\substack{i=\overline{1,m} \\a=\overline{1,n}}}$$
 which
can be also written as
$$\{e_1\wedge\ldots\wedge e_{i-1}\wedge
e_{m+a}\wedge e_{i+1}\wedge\ldots\wedge e_m\}_
{\substack{i=\overline{1,m} \\a=\overline{1,n}}}.
$$
The Riemannian
metric on the Grassmannian $G(m,n)$ is given by requesting that the
basis $\{e^*_i\otimes e_{m+a}\}$ is an orthonormal basis.

The curvature tensor field can be determined by identifying the
Grassmannian as a symmetric space (see, for example,
\cite[p.219]{PP}) and in our formalism
\begin{eqnarray}\label{eq: curv_Grassmann}
R_{P}(\rho_1,\rho_2)\rho_3 &=& \langle X_1,X_2\rangle
\langle\eta_2,\eta_3\rangle X_3^*\otimes\eta_1 -\langle
X_1,X_2\rangle \langle\eta_1,\eta_3\rangle
X_3^*\otimes\eta_2\nonumber\\
&+&\langle X_2,X_3\rangle \langle\eta_1,\eta_2\rangle
X_1^*\otimes\eta_3-\langle X_1,X_3\rangle
\langle\eta_1,\eta_2\rangle X_2^*\otimes\eta_3,
\end{eqnarray}
where $P\in G(m,n)$ and $\rho_i=X_i^*\otimes\eta_i$, $X_i\in P$,
$\eta_i\in P^\perp$, $i=1,2,3$.

We recall that for the pull-back bundle of the tangent bundle
induced by $\gamma$ we have the isometric identification
$$
\gamma^{-1}(TG(m,n))=\bigcup_{p\in M}T_{\gamma(p)}G(m,n)
=\bigcup_{p\in M}(T^*_pM\otimes N_pM)=T^*M\otimes NM,
$$
where $NM$ denotes the normal bundle of $M$ in $\r^{m+n}$.

\subsection{The tension field of the Gauss map}
In the following we intend to recall the fundamental technical steps
needed for the characterization of the harmonicity of the Gauss map.

Consider $v\in T_pM$. In order to compute $d\gamma(v)$ consider
$\sigma:I\to M$ to be a curve with $\sigma(0)=p$ and
$\dot{\sigma}(0)=v$. Let now $\{e_i\}_{i=1}^m$ be a positive
oriented orthonormal basis in $T_pM$. By parallel transporting it
along $\sigma$ we obtain  a positive oriented orthonormal basis
$\{e_i(t)\}_{i=1}^m$ in $T_{\sigma(t)}M$, for all $t$. Since
$e_i(t)$ are obtained by parallel transport along $\sigma$, we have
$$
\dot{e}_i(t)=\nabla^{\r^{m+n}}_{\dot{\sigma}}e_i=
\nabla^{M}_{\dot{\sigma}}e_i+B(\dot{\sigma},
e_i(t))=B(\dot{\sigma}, e_i(t)),
$$
where $B$ is the second fundamental form of $M$ in $\r^{m+n}$. This
implies that
\begin{eqnarray}\label{eq: dif_Gauss_map}
d\gamma_p(v)&=&\frac{d}{dt}\Big|_{t=0}(\gamma\circ\sigma)(t)=
\frac{d}{dt}\Big|_{t=0}(e_i(t)\wedge\ldots\wedge
e_m(t))\nonumber\\
&=&\sum_{i=1}^m e^*_i\otimes B(v,e_i),
\end{eqnarray}
by using the standard identifications.

The fundamental result concerning the harmonicity of the Gauss map
was obtained in \cite{RV}. We shall present here a computation
that follows \cite{EL2}. By using \eqref{eq: dif_Gauss_map} one
can compute the tension field of the Gauss map in terms of the
mean curvature of $M$. Since the bundles $\gamma^{-1}(TG(m,n))$
and $T^*M\otimes NM$ are isometric, we can write
$$
\nabla^\gamma \rho=\nabla\omega\otimes\eta+\omega\otimes
\nabla^\perp \eta,
$$
where the section $\rho\in C(\gamma^{-1}(TG(m,n)))$ in the pull-back
bundle is such that it can naturally be identified with
$\omega\otimes \eta\in C(T^*M\otimes NM)$.

Consider $\{E_i\}_{i=1}^m$ to be a local positive oriented
orthonormal frame field, geodesic at $p\in M$. By using the
expression \eqref{eq: dif_Gauss_map} for the differential of the
Gauss map and the consequence of the Codazzi equation,
$\nabla^\perp_{E_i} B(E_j,E_k)=\nabla^\perp_{E_j} B(E_i,E_k)$, for
all $ i, j, k$, we get at $p$,
\begin{eqnarray}\label{eq: tens_field_Gauss_map}
\tau(\gamma)&=&\sum_{i=1}^m\nabla d\gamma(E_i,E_i)=m\sum_{j=1}^m
E_j^*\otimes \nabla^\perp_{E_j} H,
\end{eqnarray}
where $H=\frac{1}{m}\trace B$ is the mean curvature vector field of
$M$ in $\r^{m+n}$. We note that $E_j^*$ coincides with $E_j^{\flat}$
obtained by the musical isomorphism $\flat$.

Then, the Ruh-Vilms Theorem is an immediate consequence of
~\eqref{eq: tens_field_Gauss_map}.

\section{The biharmonic equation for the Gauss map} \qquad

Inspired by the expression for the tension field given in the
previous section, we now characterize the biharmonicity of the Gauss
map in terms of the second fundamental form of the submanifold. We
obtain
\begin{theorem}\label{th: bih_gauss}
The Gauss map associated to a $m$-dimensional orientable submanifold
$M$ of $\r^{m+n}$ is biharmonic if and only if
\begin{eqnarray}\label{eq: caract_bih_Gauss_map}
\nabla^\perp_X\Delta^\perp H-m\nabla^\perp_{A_H(X)}H +\trace
B\big( 2A_{\nabla^\perp_{(\,\cdot\,)}H}(X)-A_
{\nabla^\perp_{X}H}(\,\cdot\,),\,\cdot\,\big)&\nonumber\\
-2\trace R^\perp(\,\cdot\,,X)\nabla^\perp_{\,\cdot\,}H
-\trace(\nabla^\perp_{\,\cdot\,}R^\perp)(\,\cdot\,,X)H&=&0,
\end{eqnarray}
for all $X\in C(TM)$, where $A$ denotes the Weingarten operator
and $H$ the mean curvature vector field of $M$ in $\r^{m+n}$.
\end{theorem}

\begin{proof}
We fix an orientation on $M$ and consider $\{E_i\}_{i=1}^m$ to be a
local positive oriented orthonormal frame field, geodesic at $p\in
M$. In order to determine the bitension field of the Gauss map, by
using \eqref{eq: curv_Grassmann}, \eqref{eq: dif_Gauss_map} and
\eqref{eq: tens_field_Gauss_map}  we obtain
\begin{eqnarray*}
\trace R(d\gamma,\tau(\gamma))d\gamma&=&\sum_{i=1}^m
R(d\gamma(E_i),\tau(\gamma))d\gamma(E_i)\nonumber\\
&=&m\sum_{h,i,j,k=1}^mR\big(E_j^*\otimes B(E_i,E_j),E_h^*\otimes
\nabla^\perp_{E_h}H\big) E_k^*\otimes B(E_i,E_k)\nonumber
\\
&=&m\sum_{h,i,j,k=1}^m\Big\{\delta_{jh}
\big\langle\nabla^\perp_{E_h}H,B(E_i,E_k)\big\rangle
E_k^*\otimes B(E_i,E_j)\nonumber\\
&&-\delta_{jh}\big\langle B(E_i,E_j),B(E_i,E_k)\big\rangle
E_k^*\otimes \nabla^\perp_{E_h}H\nonumber\\
&&+\delta_{hk}\big\langle\nabla^\perp_{E_h}H,B(E_i,E_j)\big\rangle
E_j^*\otimes B(E_i,E_k)\nonumber\\
&&-\delta_{jk}\big\langle\nabla^\perp_{E_h}H,B(E_i,E_j)\big\rangle
E_h^*\otimes B(E_i,E_k)\nonumber\Big\}\\
&=&m\sum_{i,j,k=1}^m\Big\{\
\big\langle\nabla^\perp_{E_j}H,B(E_i,E_k)\big\rangle
E_k^*\otimes B(E_i,E_j)\nonumber\\
&&-\big\langle B(E_i,E_j),B(E_i,E_k)\big\rangle
E_k^*\otimes \nabla^\perp_{E_j}H\nonumber\\
&&+\big\langle\nabla^\perp_{E_k}H,B(E_i,E_j)\big\rangle
E_j^*\otimes B(E_i,E_k)\nonumber\\
&&-\big\langle\nabla^\perp_{E_k}H,B(E_i,E_j)\big\rangle E_k^*\otimes
B(E_i,E_j)\nonumber\Big\}.
\end{eqnarray*}
Further computations lead to
\begin{eqnarray*}
\trace R(d\gamma,\tau(\gamma))d\gamma&=&m\sum_{k=1}^m
\Big\{E_k^*\otimes \sum_{i,j=1}^m \Big(\big(2\big\langle
\nabla^\perp_{E_j}H,B(E_i,E_k)\big\rangle\nonumber\\&&-
\big\langle B(E_i,E_j),\nabla^\perp_{E_k}H
\big\rangle\big)B(E_i,E_j)\nonumber\\
&&-\big\langle B(E_i,E_j),B(E_i,E_k)\big\rangle
\nabla^\perp_{E_j}H\Big)\Big\}.
\end{eqnarray*}
By using the Weingarten operator we can express
\begin{eqnarray*}
2\big\langle\nabla^\perp_{E_j}H,B(E_i,E_k)\big\rangle -\big\langle
B(E_i,E_j),\nabla^\perp_{E_k}H \big\rangle=\big\langle
2A_{\nabla^\perp_{E_j}}(E_k)-A_
{\nabla^\perp_{E_k}}(E_j),E_i\big\rangle,
\end{eqnarray*}
and from the Gauss equation of $M$ in $\r^{m+n}$,
$$
\big\langle B(Y,T),B(X,Z)\big\rangle=\big\langle
B(X,T),B(Y,Z)\big\rangle-\big\langle
R^M(X,Y)Z,T\big\rangle,\,\forall X,Y,Z,T\in C(TM),
$$
for $Y=Z=E_i$, $X=E_k$, $T=E_j$, we get
\begin{eqnarray*}
\sum_{i,j=1}^m \big\langle B(E_i,E_j),B(E_i,E_k)\big\rangle
E_j&=&\sum_{i,j=1}^m\Big\{ \big\langle
B(E_k,E_j),B(E_i,E_i)\big\rangle-\big\langle
R^M(E_k,E_i)E_i,E_j\big\rangle\Big\}E_j\nonumber\\
&=&\sum_{i=1}^m\Big\{m\big\langle A_H(E_k),E_i)\big\rangle
E_i -R^M(E_k,E_i)E_i \Big\}\nonumber\\
&=&mA_H(E_k)-\ricci^M(E_k),
\end{eqnarray*}
where $\ricci^M$ denotes the Ricci tensor field of $M$. \\By summing
up all of the above we obtain
\begin{eqnarray}\label{eq: bitens_Gauss_map (1)}
\trace R(d\gamma,\tau(\gamma))d\gamma&=&m\sum_{k=1}^m
E_k^*\otimes\Big(\sum_{j=1}^m B\big(
2A_{\nabla^\perp_{E_j}H}(E_k)-A_
{\nabla^\perp_{E_k}H}(E_j),E_j\big)\nonumber\\
&&-m\nabla^\perp_{A_H(E_k)}H+\nabla^\perp_{\ricci^M(E_k)}H\Big).
\end{eqnarray}
In order to compute $-\Delta^\gamma\tau(\gamma)$ we recall that,
since $\{E_i\}_{i=1}^m$ is geodesic at $p$,
$(\nabla_{E_i}E_k)_p=0$ and $(\nabla_{E_i}\nabla_{E_i}E_k)_p=0$,
for all $i,k=1,\ldots, m$. Thus, at $p$ we have
\begin{eqnarray}\label{eq: delta_tau_Gauss}
\trace\nabla^2\tau(\gamma)&=&\sum_{i=1}^m\nabla^\gamma_{E_i}
\nabla^\gamma_{E_i}\tau(\gamma)=m\sum_{i,k=1}^m\nabla^\gamma_{E_i}
\nabla^\gamma_{E_i}(E_k^*\otimes\nabla^\perp_{E_k}H)\nonumber\\
&=&m\sum_{i,k=1}^m\nabla^\gamma_{E_i}
\big(\nabla_{E_i}E_k^*\otimes\nabla^\perp_{E_k}H
+E_k^*\otimes\nabla_{E_i}^\perp\nabla^\perp_{E_k}H\big)\nonumber\\
&=&m\sum_{i,k=1}^m
\big(\nabla_{E_i}\nabla_{E_i}E_k^*\otimes\nabla^\perp_{E_k}H
+E_k^*\otimes\nabla^\perp_{E_i}
\nabla_{E_i}^\perp\nabla^\perp_{E_k}H\big)\nonumber\\
&=&m\sum_{i,k=1}^m\Big((\nabla_{E_i}\nabla_{E_i}E_k)^*\otimes\nabla^\perp_{E_k}H
+E_k^*\otimes\nabla^\perp_{E_i}
\nabla_{E_i}^\perp\nabla^\perp_{E_k}H\Big)\nonumber\\
&=&m\sum_{i,k=1}^m E_k^*\otimes\nabla^\perp_{E_i}
\nabla_{E_i}^\perp\nabla^\perp_{E_k}H.
\end{eqnarray}
Moreover, by using the curvature tensor fields $R^M$ and
$R^\perp$, for $\nabla$ and $\nabla^\perp$, respectively, at $p$
we obtain
\begin{eqnarray}\label{eq: bitens_Gauss_map (3)}
\nabla^\perp_{E_i}\nabla_{E_i}^\perp\nabla^\perp_{E_k}H&=&
\nabla^\perp_{E_i}\left(R^\perp(E_i,E_k)H+
\nabla^\perp_{E_k}\nabla^\perp_{E_i}H+
\nabla^\perp_{[E_i,E_k]}H\right)\nonumber\\
&=&(\nabla^\perp_{E_i}R^\perp)(E_i,E_k)H
+R^\perp(E_i,E_k)\nabla^\perp_{E_i}H\nonumber\\&&+
\nabla^\perp_{E_i}\nabla^\perp_{E_k}\nabla^\perp_{E_i}H+
\nabla^\perp_{E_i}\nabla^\perp_{[E_i,E_k]}H\nonumber\\
&=&(\nabla^\perp_{E_i}R^\perp)(E_i,E_k)H
+2R^\perp(E_i,E_k)\nabla^\perp_{E_i}H\nonumber\\
&&+\nabla^\perp_{E_k}\left(\nabla^\perp_{E_i}\nabla^\perp_{E_i}H
-\nabla^\perp_{\nabla_{E_i}E_i}H\right)+\nabla^\perp_{R^M(E_k,E_i)E_i}H\nonumber,
\end{eqnarray}
thus
\begin{eqnarray}\label{eq: bitens_Gauss_map (2)}
\trace\nabla^2\tau(\gamma)&=&m\sum_{k=1}^m
E_k^*\otimes\Big(\trace\big\{(\nabla^\perp_{\,\cdot\,}R^\perp)(\,\cdot\,,E_k)H
+2R^\perp(\,\cdot\,,E_k)\nabla^\perp_{\,\cdot\,}H\big\}\nonumber
\\&&-\nabla^\perp_{E_k}\Delta^\perp
H+\nabla^\perp_{\ricci^M(E_k)}H\Big).
\end{eqnarray}
We finally substitute \eqref{eq: bitens_Gauss_map (1)} and
\eqref{eq: bitens_Gauss_map (2)} into the biharmonic equation and
conclude.
\end{proof}

\subsection{The case of hypersurfaces}

Let $M$ be a nowhere zero mean curvature hypersurface in $\r^{m+1}$.
We obtain

\begin{theorem}\label{th: bih_Gauss_hypersurface}
The Gauss map of a non-minimal hypersurface $M^m$ in $\r^{m+1}$ is
proper biharmonic if and only if $\grad f\neq 0$ and
\begin{equation}\label{eq: caract_bih_Gauss_map_hyper}
\Delta\grad f +A^2(\grad f)-|A|^2\grad f=0.
\end{equation}
where $\Delta$ denotes the rough Laplacian on $C(TM)$ and $f$ and
$A$ denote the mean curvature function and, respectively, the shape
operator of $M$ in $\r^{m+1}$.
\end{theorem}

\begin{proof}
In this case we can consider the expression of the mean curvature
vector field as $H=f\eta$, where $f=\vert H\vert$ is the mean
curvature function of $M$ in $\r^{m+1}$ and $\eta=\frac{1}{\vert
H\vert}H$ is a unit section in the normal bundle $NM$.

In the following we shall use the general biharmonic equation
\eqref{eq: caract_bih_Gauss_map} and express it for the case of
hypersurfaces.

Since $H=f\eta$ and $\nabla^\perp\eta=0$, we get that
\begin{equation}\label{eq: Gauss_hyper_1}
\trace R^\perp(\,\cdot\,,X)\nabla^\perp_{\,\cdot\,}H=0
\end{equation}
and
\begin{equation}\label{eq: Gauss_hyper_2}
\trace(\nabla^\perp_{\,\cdot\,}R^\perp)(\,\cdot\,,X)H=0.
\end{equation}
By considering now $\{E_i\}_{i=1}^m$ to be a local orthonormal frame
field on $M$ we obtain
\begin{eqnarray}\label{eq: Gauss_hyper_3}
\trace B\big( 2A_{\nabla^\perp_{(\,\cdot\,)}H}(X)-A_
{\nabla^\perp_{X}H}(\,\cdot\,),\,\cdot\,\big)&
=&\sum_{i=1}^m\Big\{2B\big(A(X),E_i(f)E_i\big)-X(f)B\big(A(E_i),E_i\big)\Big\}\nonumber\\
&=&\langle X, 2A^2(\grad f)-|A|^2\grad f\rangle\eta
\end{eqnarray}
and
\begin{equation}\label{eq: Gauss_hyper_4}
\nabla^\perp_{A_H(X)}H=\langle X, fA(\grad f)\rangle\eta.
\end{equation}
For the final term, we recall that $\grad \Delta f=\Delta\grad
f+\ricci^M(\grad f)$. Also, as a consequence of the Gauss equation
for $M$ in $\r^{m+1}$, we have
$$\ricci^M(\grad f)=mfA(\grad
f)-A^2(\grad f).
$$
Thus, it follows that
\begin{eqnarray}\label{eq: Gauss_hyper_5}
\nabla^\perp_X\Delta^\perp H&=&\langle X,\grad\Delta
f\rangle\eta\nonumber\\
&=&\langle X,\Delta\grad f+m f A(\grad f)-A^2(\grad f)\rangle\eta.
\end{eqnarray}
Finally, by replacing the expressions \eqref{eq: Gauss_hyper_1},
\eqref{eq: Gauss_hyper_2}, \eqref{eq: Gauss_hyper_3}, \eqref{eq:
Gauss_hyper_4}, \eqref{eq: Gauss_hyper_5} in \eqref{eq:
caract_bih_Gauss_map}, we deduce that the Gauss map of the
hypersurface is biharmonic if and only if
$$
\langle X, \Delta\grad f +A^2(\grad f)-|A|^2\grad
f\rangle=0,\qquad \forall X\in C(TM),
$$
and this completes the proof.
\end{proof}

\begin{example}
The first example of hypersurface with proper biharmonic Gauss map
is obtained by analyzing the right cylinders in $\r^3$. A right
cylinder in $\r^3$ determined by a plane curve parametrized by arc
length has proper biharmonic Gauss map if and only if the curve is a
{\it clothoid} or a {\it Cornu spiral}. Indeed, consider a right
cylinder in $\r^3$, determined by a curve $\sigma:I\to\r^2$,
parametrized by arc length. Denote by $s$ the parameter along
$\sigma$ and by $t$ the parameter along the generatrix. If $k$
denotes the signed curvature of $\sigma$, then the mean curvature
function is $f(s,t)=\pm\frac{1}{2}k(s)$. Moreover, since $A^2(\grad
f)=|A|^2\grad f=k^2(s)\grad f$, equation \eqref{eq:
caract_bih_Gauss_map_hyper} becomes $\overset{...}{k}=0$, hence $k$
is a second degree polynomial in $s$ and, by using the fundamental
theorem of plane curves, we deduce that $\sigma$ is a clothoid.
\end{example}

By a straightforward computation we obtain

\begin{proposition}
Let $M^m$ be a submanifold in $\r^{m+n}$. Then the generalized
cylinder $\r\times M$ in $\r^{m+n+1}$ has biharmonic Gauss map if
and only if $M$ has biharmonic Gauss map.
\end{proposition}

\section{Hypercones with  biharmonic Gauss map}

In order to obtain some examples of hypersurfaces with biharmonic
Gauss map, in the following we shall study the hypercones
generated by hypersurfaces of the unit Euclidean sphere.

Let us first consider $\overline{M}$ to be an arbitrary
$r$-dimensional submanifold of the unit Euclidean sphere $\s^{m+1}$.
The cone in $\r^{m+2}$ generated by $\overline{M}$ is defined by the
immersion
\begin{align*}
\F:(0,\infty)\times \overline{M}&\longrightarrow \r^{m+2}\\
(t,p)&\longmapsto t\cdot p.
\end{align*}
The differential of $\F$ is determined by
$$
d\F_{(t,p)}(\t_{(t,p)})=p= x^\alpha(p)e_\alpha(\F(p)),
$$
$$
d\F_{(t,p)}(X_{(t,p)})=tX_ p=t\xi^\alpha(p)e_\alpha(\F(p)),
$$
where $\{e_\alpha\}_{\alpha=1}^{m+2}$ is the canonical orthonormal
frame field on $\r^{m+2}$, $X\in C(T\overline{M})$ with
$X(p)=\xi^\alpha(p)e_\alpha(p)\in \r^{m+2}$, for all $p\in
\overline{M}$, and, typically, we use the same notation for a
vector field and for its lift to the product manifold.

If we denote by $\overline{g}$ the metric on $\overline{M}$, then
the immersion $\F:(0,\infty)\times \overline{M}\to\r^{m+2}$ induces
on the product $(0,\infty)\times \overline{M}$ the warped metric
$g=dt^2+t^2\overline{g}$. Thus the cone can be seen as the warped
product
$$
M=(0,\infty)\times_{t^2}\overline{M}.
$$

Denote by $\nabla$ and $\overline{\nabla}$ the Levi-Civita
connections on $M$ and on $\overline{M}$, respectively, and recall
that (see \cite[p.206]{BON}) $\nabla$ is completely determined by
\begin{equation}\label{eq: nabla_cone}
\left\{%
\begin{array}{ll}
     \nabla_{\tt}\t=0\\
     \nabla_{\tt} X=\nabla_X \t=\frac{1}{t} \, X\\
     \nabla_X Y=\overline{\nabla}_X Y-t\langle X, Y\rangle \t,
\end{array}%
\right.
\end{equation}
where $X,Y\in C(T\overline{M})$. The second fundamental form of
the cone in $\r^{m+2}$, obtained by using \eqref{eq: nabla_cone},
is given by
\begin{eqnarray}\label{eq: a doua forma fund}
B(\t,\t)=0, \qquad B(X,\t)=0,\qquad B(X,Y)=t\overline{B}(X,Y),
\end{eqnarray}
for all $X,Y\in C(T\overline{M})$, where $\overline{B}$ denotes
the second fundamental form of $\overline{M}$ in $\s^{m+1}$.

Thus, if we denote by $\overline{A}$ the Weingarten operator of
$\overline{M}$ in $\s^{m+1}$ with respect to an arbitrary fixed
unit section $\overline{\eta}$ in the normal bundle of
$\overline{M}$ in $\s^{m+1}$, we obtain the expression for the
Weingarten operator $A$ of the cone with respect to the unit
section $\eta(t,p)=\overline{\eta}(p)$, $(t,p)\in M$, in the
normal bundle of $M$ in $\r^{m+2}$,
\begin{equation}\label{eq: shape_op_hypercone}
A(\t)=0\qquad\text{and}\qquad
A(X)=\frac{1}{t}\overline{A}(X),\qquad
\end{equation}
for all $X\in C(T\overline{M})$, and, consequently, $
|A|^2=\frac{1}{t^2}|\overline{A}|^2$.

Moreover, for a smooth function $f\in C^\infty(M)$, we have
\begin{equation}\label{eq:grad_cone}
\grad f=\frac{\partial f}{\partial t}\t+\frac{1}{t^2}\grad f_t
\end{equation}
and
\begin{equation}\label{eq:Laplace_cone}
\Delta f=-\frac{\partial^2 f}{\partial
t^2}-\frac{n}{t}\frac{\partial f}{\partial
t}+\frac{1}{t^2}\overline{\Delta} f_t,
\end{equation}
where $f_t\in C^\infty(\overline{M})$, $f_t(p)=f(t,p)$, for all
$p\in \overline{M}$ and $t\in(0,\infty)$.

We are now ready to write down the conditions for the biharmonicity
of the Gauss map associated to a hypercone.

\begin{theorem}\label{th: bih_Gauss_hypercone}
Let $\overline{M}$ be a non-minimal hypersurface of $\s^{m+1}$. The
Gauss map associated to the hypercone $(0,\infty)\times_{t^2}
\overline{M}$ is proper biharmonic if and only if
\begin{equation}\label{eq: cond_gaussmap_bih_general}
\left\{
\begin{array}{l}
\overline{\Delta}\grad \overline{f}+\overline{A}^2(\grad
\overline{f})+(2m-3-|\overline{A}|^2)\grad \overline{f}=0
\\ \mbox{} \\
3\overline{\Delta}\,
\overline{f}+(3m-6-|\overline{A}|^2)\overline{f}=0,
\end{array}
\right.
\end{equation}
where $\overline{A}$ and  $\overline{f}\in C^\infty(\overline{M})$
are the shape operator and the mean curvature function of
$\overline{M}$ in $\s^{m+1}$, respectively.
\end{theorem}

\begin{proof} Consider $\t$ and $\{E_i\}_{i=1}^m$ a local orthonormal frame field on
$\overline{M}$, geodesic at $p$. Then
$\Big\{\t,\frac{1}{t}E_i\Big\}_{i=1}^m$ constitutes a local
orthonormal frame field on $(0,\infty)\times_{t^2}\overline{M}$.
Denoting by $\overline{f}$ is the mean curvature function of
$\overline{M}$ in $\s^{m+1}$ and using \eqref{eq: a doua forma
fund}, we get the mean curvature function $f$ of the hypercone,
$$
f=\frac{m}{(m+1)}\frac{1}{t}\overline{f}.
$$
Using \eqref{eq:grad_cone}, we get
$$
\grad f=\frac{m}{(m+1)}\Big(-\frac{1}{t^2}\overline{f}\t+
\frac{1}{t^3}\grad \overline{f}\Big),
$$
and this, together with \eqref{eq: shape_op_hypercone}, implies
\begin{eqnarray}\label{eq: bitens_part1}
A^2(\grad f)=\frac{m}{(m+1)} \frac{1}{t^5}\overline{A}^2(\grad
\overline{f}).
\end{eqnarray}
Also,
\begin{equation}\label{eq: bitens_part2}
-|A|^2\grad f=\frac{m}{(m+1)}
\frac{1}{t^4}|\overline{A}|^2\Big(\overline{f}\t-\frac{1}{t}\grad
\overline{f}\Big).
\end{equation}
In order to compute $\Delta(\grad f)=-\trace \nabla^2(\grad f)$ we
shall use \eqref{eq: nabla_cone}. Thus,
\begin{eqnarray*}
\trace\nabla^2\Big(\frac{1}{t^2}\overline{f}\t\Big)
&=&\overline{f}\nabla_{\tt} \nabla_{\tt}\Big(\frac{1}{t^2}\t\Big)
+\frac{1}{t^2}\sum_{i=1}^m\Big\{\frac{1}{t^2}\nabla_{E_i}\nabla_{E_i}
(\overline{f} \t)\\&&-\nabla_{\nabla_{E_i}E_i}
\Big(\frac{1}{t^2}\overline{f}\t\Big)\Big\}\nonumber\\
&=&\frac{6}{t^4}\overline{f}\t + \frac{1}{t^4}\sum_{i=1}^m\Big\{
E_i(E_i(\overline{f}))\t+\frac{2}{t}E_i(\overline{f})E_i-3\overline{f}\t\Big\}\nonumber\\
&=&\frac{1}{t^4}\big((6-3m)\overline{f}-\overline{\Delta}\,
\overline{f}\big)\t+\frac{2}{t^5}\grad \overline{f}\nonumber,
\end{eqnarray*}
and
\begin{eqnarray*}
\trace\nabla^2\Big(\frac{1}{t^3}\grad \overline{f}\Big)
&=&\nabla_{\tt} \nabla_{\tt}\Big(\frac{1}{t^3}\grad
\overline{f}\Big)
\\&&+\frac{1}{t^2}\sum_{i=1}^m\Big\{\frac{1}{t^3}
\nabla_{E_i}\nabla_{E_i} \grad
\overline{f}-\nabla_{\nabla_{E_i}E_i}
\Big(\frac{1}{t^3}\grad \overline{f}\Big)\Big\}\\
&=&\frac{6}{t^5}\grad \overline{f}
+\frac{1}{t^2}\sum_{i=1}^m\Big\{\frac{1}{t^3}
\nabla_{E_i}\big(\overline{\nabla}_{E_i}\grad \overline{f}-
t\langle E_i,\grad \overline{f}\rangle\t\big)\\
&&-\frac{2}{t^3}\grad \overline{f}\Big\}\\
&=&\frac{2}{t^4}\overline{\Delta}\,
\overline{f}\t+\frac{1}{t^5}\big((5-2m)\grad
\overline{f}-\overline{\Delta}\grad \overline{f}\big).
\end{eqnarray*}
Using the two expressions above we obtain
\begin{equation}\label{eq: bitens_part3} \Delta(\grad f)=
\frac{m}{m+1}\Big\{\frac{1}{t^4}\big((6-3m)\overline{f}-3\overline{\Delta}\,
\overline{f}\big)\t+\frac{1}{t^5}\big((2m-3)\grad
\overline{f}+\overline{\Delta}\grad \overline{f}\big)\Big\}.
\end{equation}
By substituting \eqref{eq: bitens_part1}, \eqref{eq: bitens_part2}
and \eqref{eq: bitens_part3} in \eqref{eq:
caract_bih_Gauss_map_hyper}, we obtain the desired result.
\end{proof}

\begin{corollary}
Let $\overline{M}$ be a non-minimal hypersurface of $\s^{m+1}$ with
constant norm of the shape operator. We have
\begin{itemize}
\item [(i)] if the Gauss map associated to the hypercone
$(0,\infty)\times_{t^2} \overline{M}$ is proper biharmonic, then
\begin{equation}\label{eq: cor1_2}
2\overline{A}^2(\grad \overline{f})-m\overline{f}\,\overline{A}(\grad
\overline{f})-\frac{2}{3}|\overline{A}|^2\grad \overline{f}=0.
\end{equation}
\item [(ii)]if $\overline{M}$ is compact, then the Gauss map
associated to the hypercone is proper biharmonic if and only if
$\overline{M}$ has constant mean curvature in $\s^{m+1}$, $m>2$ and
$\vert\overline{A}\vert^2=3(m-2)$.
\end{itemize}
\end{corollary}

\begin{proof}
The second equation of \eqref{eq: cond_gaussmap_bih_general}
implies
$$
3\grad \overline{\Delta}\,
\overline{f}+(3m-6-|\overline{A}|^2)\grad \overline{f}=0,
$$
and, since $\grad \overline{\Delta}\,
\overline{f}=\overline{\Delta}\grad
\overline{f}+\ricci^{\overline{M}}(\grad \overline{f})$, we obtain
$$
\overline{\Delta}\grad
\overline{f}=\left(2-m+\frac{1}{3}|\overline{A}|^2\right)\grad\overline{f}-\ricci^{\overline{M}}(\grad
\overline{f}).
$$
We substitute this expression in the first equation of \eqref{eq:
cond_gaussmap_bih_general} and it follows that
$$
\overline{A}^2(\grad\overline{f})-\ricci^{\overline{M}}(\grad
\overline{f})+\left(m-1-\frac{2}{3}|\overline{A}|^2\right)\grad
\overline{f}=0.
$$
Finally, from the Gauss equation for $\overline{M}$ in $\s^{m+1}$
we deduce that
$$
\ricci^{\overline{M}}(X)=(m-1)X+m\overline{f}\,\overline{A}(X)-\overline{A}^2(X),\qquad
\forall X\in C(T\overline{M}),
$$
and we conclude.

In order to prove (ii), we integrate the second equation of
\eqref{eq: cond_gaussmap_bih_general} and we get
$3m-6-\vert\overline{A}\vert^2=0$, and then $\overline{f}$ is
constant.
\end{proof}

As a consequence of Theorem \ref{th: bih_Gauss_hypercone} we
obtain,
\begin{theorem}\label{th: CMC_hypercone}
Let $\overline{M}$ be a constant non-zero mean curvature
hypersurface of $\s^{m+1}$. The Gauss map associated to the
hypercone $(0,\infty)\times_{t^2} \overline{M}$ is proper biharmonic
if and only if $m>2$ and $|\overline{A}|^2=3(m-2)$, where
$\overline{A}$ is the shape operator of $\overline{M}$ in
$\s^{m+1}$.
\end{theorem}

\begin{proof}
If $\overline{f}$ is constant, then the first condition of
\eqref{eq: cond_gaussmap_bih_general} is identically satisfied and
the second one implies $|\overline{A}|^2=3(m-2)$. The converse is
immediate.
\end{proof}

For the case of hypercones in $\r^3$ and $\r^4$, i.e. $m=1$ and
$m=2$, we have the following non-existence results.

\begin{theorem}\label{th: non-ex_cone_Gauss_R3}
There exist no cones in $\r^3$ with proper biharmonic Gauss map.
\end{theorem}

\begin{proof}
Consider a cone in $\r^3$ generated by a curve $\sigma:I\to\s^2$,
parametrized by arc length. Denote by $s$ the parameter on the curve
and by $T=\dot{\sigma}$ the tangent vector field along $\sigma$.
Since $\nabla^{\ts^2}_T T=kN$, with $N$ the unit normal vector field
along $\sigma$, the mean curvature function is given by
$\overline{f}=\pm k$,
$$
\overline{A}(\ss)=k\ss\qquad\text{and}\qquad |\overline{A}|^2=k^2,
$$
$$
\grad \overline{f}=\pm\dot{k}\ss\qquad\text{and}\qquad
\overline{\Delta}\grad \overline{f}=\mp\overset{...}{k}\ss.
$$
Thus, condition \eqref{eq: cond_gaussmap_bih_general} becomes
\begin{equation}\label{eq: cond_gaussmap_bih_conR3}
\left\{
\begin{array}{l}
\overset{...}{k}+\dot{k}=0
\\ \mbox{} \\
k(3+k^2)+3\ddot{k}=0.
\end{array}
\right.
\end{equation}
This implies that $\dot{k}k^2=0$, hence $k=0$, i.e. the Gauss map of
the cone is harmonic, and we conclude.
\end{proof}

\begin{theorem}\label{th: non-ex_cone_Gauss_R4}
There exist no hypercones in $\r^4$, over compact non-minimal
surfaces $\overline{M}^2\subset\s^3$, with proper biharmonic Gauss
map.
\end{theorem}

\begin{proof}
Suppose that the Gauss map of the hypercone over $\overline{M}$ is
proper biharmonic. Since $m=2$, the second equation of \eqref{eq:
cond_gaussmap_bih_general} leads to
\begin{equation}\label{eq: consec_hyper_R4}
3\overline{\Delta}\,\overline{f}-|\overline{A}|^2\overline{f}=0.
\end{equation}
By integrating condition \eqref{eq: consec_hyper_R4} on
$\overline{M}$ and by using the fact that $\overline{f}$ is
positive, we conclude that $|\overline{A}|^2=0$ and we have a
contradiction.
\end{proof}

When $m>2$, we have examples of hypercones with proper biharmonic
Gauss map. Recall that if $\overline{M}$ is a hypersurface in
$\s^{m+1}$, then the cone over $\overline{M}$ has harmonic Gauss map
if and only if $\overline{M}$ is minimal in $\s^{m+1}$ (see
\cite{OJG,JS}). This does not hold in the case of the biharmonicity.
Indeed, by considering $\overline{M}$ to be a constant mean
curvature proper biharmonic hypersurface of $\s^{m+1}$ and by using
the fact that the squared norm of the shape operator of such a
submanifold is equal to $m$ (see \cite{RCSMCO1}) we get

\begin{theorem}\label{cor:hypercone_bih_hypersurf_CMC}
Let $\overline{M}$ be a constant mean curvature proper
biharmonic hypersurface of $\s^{m+1}$. Then the hypercone
$(0,\infty)\times_{t^2} \overline{M}$ has proper biharmonic
associated Gauss map if and only if $m=3$.
\end{theorem}

\begin{remark}
Theorem~\ref{cor:hypercone_bih_hypersurf_CMC} can be deduced, for
the particular case of the hypersphere of radius equal to $\rrm$, in
a more geometrical manner. The argument is the following. In
\cite{LO2}, the authors proved that if $\psi:N\to\s^{m}(\rrm)$ is a
harmonic map and ${\bf i}:\s^{m}(\rrm)\to\s^{m+1}$ denotes the
inclusion map, then the tension and bitension fields of the
composition are given by
\begin{equation}\label{eq: expr_bitension_comp}
\tau({\bf i}\circ\psi)=-2e(\psi)\overline{\eta}
\quad \hbox{and}
\quad \frac{1}{2}\tau_2({\bf i}\circ\psi)=\big(\Delta
e(\psi)\big)\overline{\eta}-2d\psi(\grad e(\psi)),
\end{equation}
where $e(\psi)$ denotes the energy density of the map $\psi$ and
$\overline{\eta}$ the unit section of the normal bundle of
$\s^m(\rrm)$ in $\s^{m+1}$. The Gauss map associated to the
hypercone $(0,\infty)\times_{t^2} \s^m(\rrm)$ are given by
$$
\gamma:(0,\infty)\times_{t^2} \s^m(\rrm)\to \s^{m+1}
$$
$$
\gamma(t,p)=\overline{\eta}(p),
$$
i.e. $\gamma(t,x^1,\ldots,x^{m+1},\rrm)=(x^1,\ldots,x^{m+1},-\rrm)$.
We can thus think of $\gamma$, up to an isometry, as the composition
${\bf i}\circ\psi$, where
$$
\psi:(0,\infty)\times_{t^2} \s^m(\rrm)\to \s^{m}(\rrm)
$$
$$
\psi(t,p)=p.
$$
The map $\psi$ is the projection onto the second factor of a
warped product, so it is a harmonic map.
Now, since $d\psi(\t)=0$ and $d\psi(X)=X$, for all $X\in
C(T\s^m(\rrm))$, the energy density of $\psi$ is $
e(\psi)=\frac{m}{2t^2}$ and $\grad e(\psi)=-\frac{m}{t^3}\t$.
Thus, we deduce that
\begin{equation}\label{eq: eqaux}
\Delta e(\psi)=\frac{m(m-3)}{t^4}\qquad\text{and}\qquad
d\psi(\grad e(\psi))=0.
\end{equation}
Finally, by using \eqref{eq: expr_bitension_comp} and \eqref{eq:
eqaux}, we conclude that the Gauss map associated to the hypercone
$(0,\infty)\times_{t^2}\s^m(\rrm)$ in $\r^{m+2}$ is proper
biharmonic if and only if $m=3$, in accordance with Theorem
\ref{cor:hypercone_bih_hypersurf_CMC}.
\end{remark}

\section{Hypercones generated by isoparametric
hypersurfaces in spheres}

We recall that a hypersurface $\overline{M}^m$ in $\s^{m+1}$ is
said to be {\it isoparametric} of type $\ell$ if it has constant
principal curvatures $k_1 > \ldots
> k_\ell$ with respective constant multiplicities $m_1, \ldots
,m_\ell$, $m = m_1 + m_2 + \ldots + m_\ell$. E.~Cartan classified
in \cite{CAR} the isoparametric hypersurfaces with $\ell=1,2,3$.
For $\ell>3$ a full classification of isoparametric hypersurfaces
is not yet known. Nevertheless, it is known that the number $\ell$
is either $1, 2, 3, 4$ or $6$ (see \cite{MUN}) and the following
information on the principal curvatures and their multiplicities
is available.
\begin{itemize}
\item[(i)] If $\ell  = 1$, then $\overline{M}$ is totally umbilical.
\item[(ii)] If $\ell = 2$, then $\overline{M} = \s^{m_1}(r_1)\times \s^{m_2}(r_2)$,
$r_1^2 + r_2^2 = 1$.
\item[(iii)] If $\ell = 3$, then $m_1 = m_2 = m_3 = 2^q$, $q=0,1,2,3$.
\item[(iv)] If $\ell = 4$, then $m_1 = m_3$ and $m_2 = m_4$. Moreover,
$(m_1,m_2)=(2,2)$ or $(4,5)$, or $m_1+m_2+1$ is a multiple of
$2^{\rho(m^{\ast}-1)}$, where $\rho(s)$ is the number of integers
$r$ with $1\leq r\leq s$, $r\equiv 0, 1, 2, 4\, (\textrm{mod}\, 8)$
and $m^{\ast}={\rm{min}}\{m_1,m_2\}$.
\item[(v)] If $\ell = 6$, then $m_1=m_2=\ldots=m_6=1$ or $2$.
\end{itemize}
Moreover, there exists an angle $\theta$, $0 < \theta <
\frac{\pi}{\ell}$ , such that
\begin{equation}\label{eq-kalpha}
k_\alpha = \cot\big(\theta + \frac{(\alpha-1)\pi}{\ell}\big),\quad
\alpha = 1,\ldots, \ell.
\end{equation}

We now study the biharmonicity of the Gauss map of the  hypercones
generated by isoparametric hypersurfaces in spheres. We shall detail
this study according to the type $\ell$ of the isoparametric
hypersurface.

\subsection*{Isoparametric hypersurface with ${\mathbf \ell=1}$}
In this case $\overline{M}$ is a hypersphere
$\s^m(a)$, $a\in(0,1)$, in $\s^{m+1}$. Since
$\displaystyle{|\overline{A}|^2=m\frac{1-a^2}{a^2}}$, by using
Theorem~\ref{th: CMC_hypercone}, we obtain

\begin{proposition}\label{cor:hypercone_hypersphere}
Consider the hypercone $(0,\infty)\times_{t^2} \s^m(a)$,
$a\in(0,1)$, in $\r^{m+2}$. Its associated Gauss map is proper
biharmonic if and only if $m>2$ and $a=\sqrt{\frac{m}{4m-6}}$.
\end{proposition}

\begin{remark}
We underline the fact that
Proposition~\ref{cor:hypercone_hypersphere} provides examples of
hypersurfaces with proper biharmonic associated Gauss map in any
$(m+2)$-dimensional Euclidean space, with $m>2$.
\end{remark}

\subsection*{Isoparametric hypersurface with $\ell=2$}
In this case $\overline{M}$ is a generalized
torus $\s^{m_1}(r_1)\times\s^{m_2}(r_2)\subset \s^{m+1}$,
$m_1+m_2=m$, $r_1^2+r_2^2=1$. The squared norm of the shape
operator is
$|\overline{A}|^2=\big(\frac{r_2}{r_1}\big)^2m_1+\big(\frac{r_1}{r_2}\big)^2m_2$
and, by using Theorem~\ref{th: CMC_hypercone}, we get

\begin{proposition}\label{cor:hypercone_CLiff_torus}
Consider the hypercone $(0,\infty)\times_{t^2}
\big(\s^{m_1}(r_1)\times\s^{m_2}(r_2)\big)$ in $\r^{m+2}$. Its
associated Gauss map is proper biharmonic if and only if $m>3$,
$\frac{m_1}{r_1^2}\neq\frac{m_2}{r_2^2}$ and
\begin{equation}\label{eq: Gauss_bih_Clifford_torus}
\frac{m_1}{r_1^2}+\frac{m_2}{r_2^2}=4m-6.
\end{equation}
\end{proposition}

\begin{remark}
In order to obtain an example, consider $m>3$, $m_1=1$, and
$m_2=m-1$. Then
$$
r_1^2=\frac{3m-4\pm \sqrt{9m^2-40m+40}}{2(4m-6)}\quad \textrm{and}
\quad r_2^2=\frac{5m-8\mp \sqrt{9m^2-40m+40}}{2(4m-6)}$$ are
solutions for \eqref{eq: Gauss_bih_Clifford_torus}.
\end{remark}

\subsection*{Isoparametric hypersurface with $\ell=3$}
In this case, taking into account
\eqref{eq-kalpha}, there exists $\theta\in(0,\pi/3)$ such that
$$
k_1=\cot \theta,\qquad k_2 = \cot\big(\theta +
\frac{\pi}{3}\big)=\frac{k_1-\sqrt 3}{1+\sqrt{3}k_1},\qquad k_3 =
\cot\big(\theta + \frac{2\pi}{3}\big)=\frac{k_1+\sqrt
3}{1-\sqrt{3}k_1}.
$$
Thus, the square of the norm of the shape operator is
\begin{equation}\label{eq: type 3 A^2}
|\overline{A}|^2=2^q(k_1^2+k_2^2+k_3^2)
=2^q\frac{9k_1^6+45k_1^2+6}{(1-3k_1^2)^2}
\end{equation}
and $m=3\cdot2^q$, $q=0,1,2,3$.

On the other hand, from Theorem~\ref{th: CMC_hypercone}, the
hypercone generated by $\overline{M}$ has proper biharmonic Gauss
map if and only if
$$
|\overline{A}|^2=3(m-2)=3(3\cdot2^q-2).
$$
The last equation, together with \eqref{eq: type 3 A^2}, implies
that $k_1$ is a solution of
\begin{equation}\label{eq: polin_l=3}
3\cdot2^q\, x^6+(18-27\cdot 2^q)\,x^4+(-12+33\cdot
2^q)\,x^2+2-2^q=0.
\end{equation}

If $q=0$, equation \eqref{eq: polin_l=3} becomes $3
x^6-9x^4+11x^2+1=0$ and it has no real roots.

If $q=1$, equation \eqref{eq: polin_l=3} becomes $x^2(x^2-3)^2=0$,
which has one root $x=\sqrt 3$ in $(\frac{\sqrt 3}{3},\infty)$.
Notice that when $k_1=\sqrt{3}$, $\overline{M}$ is minimal.

If $q=2$, equation \eqref{eq: polin_l=3} becomes
$6x^6-45x^4+60x^2-1=0$ and it has two distinct roots in
$(\frac{\sqrt 3}{3},\infty)$, different from $\sqrt 3$.

If $q=3$, equation \eqref{eq: polin_l=3} becomes
$4x^6-33x^4+42x^2-1=0$ and it has two distinct roots in
$(\frac{\sqrt 3}{3},\infty)$, different from $\sqrt 3$.

We can conclude,

\begin{proposition}\label{cor: isoparam_l=3}
Consider an isoparametric hypersurface of type $3$ in $\s^{m+1}$.
The Gauss map of its hypercone is proper biharmonic if and only if
\begin{itemize}
\item[(i)] $q=2$ and the first principal curvature $k_1$
is one of the two roots of the equation $6x^6-45x^4+60x^2-1=0$ in
$(\frac{\sqrt 3}{3},\infty)$,

or

\item[(ii)] $q=3$ and the first principal curvature $k_1$
is one of the two roots of the equation
$4x^6-33x^4+42x^2-1=0$ in $(\frac{\sqrt 3}{3},\infty)$.
\end{itemize}
\end{proposition}

\subsection*{Isoparametric hypersurface with $\ell=4$} In this case, taking into account
\eqref{eq-kalpha}, there exists $\theta\in(0,\pi/4)$ such that
$$
k_1=\cot \theta,\qquad k_2 = \cot\big(\theta +
\frac{\pi}{4}\big)=\frac{k_1-1}{k_1+1},
$$
$$
k_3 = \cot\big(\theta + \frac{\pi}{2}\big)=-\frac{1}{k_1}, \qquad
k_4 = \cot\big(\theta + \frac{3\pi}{4}\big)=-\frac{k_1+1}{k_1-1}.
$$

The square of the norm of the shape operator is
\begin{eqnarray}\label{eq: type 4 A^2}
|\overline{A}|^2&=&m_1\left(k_1^2+\frac{1}{k_1^2}\right)+
m_2\left[\left(\frac{k_1-1}{k_1+1}\right)^2+\left(\frac{k_1+1}{k_1-1}\right)^2\right]\nonumber\\
&=&m_1\lambda+16m_2\frac{1}{\lambda}+2(m_1+m_2),
\end{eqnarray}
where $\lambda=\left(k_1-\frac{1}{k_1}\right)^2$.

Since in this case $m=2(m_1+m_2)$, from Theorem~\ref{th:
CMC_hypercone}, the hypercone generated by $\overline{M}$ has
proper biharmonic Gauss map if and only if
$$
|\overline{A}|^2=3(m-2)=6(m_1+m_2-1).
$$
The last equation, together with \eqref{eq: type 4 A^2}, implies
that $\lambda$ is a solution of
\begin{equation}\label{eq: polin_l=4}
m_1\lambda^2-(4(m_1+m_2)-6)\lambda+16m_2=0.
\end{equation}

Notice that if $(m_1,m_2)=(2,2)$ or $(4,5)$, equation \eqref{eq:
polin_l=4} has no real roots. Consequently, we obtain

\begin{proposition}\label{cor: isoparam_l=4}
Consider an isoparametric hypersurface of type $4$ in $\s^{m+1}$.
The Gauss map of its hypercone is proper biharmonic if and only if
its first principal curvature $k_1$ is given by the condition that
$\lambda=\left(k_1-\frac{1}{k_1}\right)^2$ is the positive
solution of the equation
$$
m_1\lambda^2-\big(4(m_1+m_2)-6\big)\lambda+16m_2=0,
$$
and $m_1+m_2+1$ is a multiple of $2^{\rho(m^{\ast}-1)}$, where
$\rho(s)$ is the number of integers $r$ with $1\leq r\leq s$ and
$r\equiv 0, 1, 2, 4\, (\textrm{mod}\, 8)$.
\end{proposition}

\begin{example}
In order to obtain an explicit example for this case we shall
consider from the Takagi list (see \cite{T1}) the following
homogeneous hypersurfaces with four principal curvatures,
\begin{equation}\label{eq: type C}
\overline{M}=S(U(k)\times U(2))/(T^2\times
SU(k-2))\subset\s^{2n+1},
\end{equation}
where $n=2k+1$, $n\geq 5$.

Since $m_1=n-2$ and $m_2=2$, from Proposition \ref{cor:
isoparam_l=4} we deduce that the Gauss map of the hypercone over
$\overline{M}$ is biharmonic if and only if
\begin{equation}\label{eq: l=4, eq_lambda}
(n-2)\lambda^2-(4n-6)\lambda+32=0.
\end{equation}
By denoting $\sin^2 2\theta=x\in(0,1)$ we have
$\lambda=4\frac{1-x}{x}$, and equation \eqref{eq: l=4, eq_lambda}
becomes
\begin{equation}\label{eq: con_Tip C}
(4n-3)x^2-(6n-11)x+2(n-2)=0.
\end{equation}
Since $n$ is odd, \eqref{eq: con_Tip C} has real roots if
and only if $n\geq 9$. It is easy to verify that, for $n\geq 9$,
the two real roots of \eqref{eq: con_Tip C} are in $(0,1)$ and we
obtain
\begin{equation}\label{eq:TipC1_bihGauss}
\sin^2 2\theta=\frac{6n-11\pm\sqrt{4n^2-44n+73}}{2(4n-3)}.
\end{equation}
Notice that $\overline{M}$ is minimal if and only if
$$
(n-2)\cot^4\theta-2(n+2)\cot^2\theta+(n-2)=0,
$$
thus, the Gauss map of $(0,\infty)\times_{t^2}
\overline{M}$ is harmonic if and only if
\begin{equation}\label{eq:TipC2_bihGauss}
\cot^2\theta = \frac{\sqrt n\pm \sqrt 2}{\sqrt n\mp \sqrt 2}.
\end{equation}
Since $n\geq 9$, from   \eqref{eq:TipC1_bihGauss} and
\eqref{eq:TipC2_bihGauss} we get that
$(0,\infty)\times_{t^2} \overline{M}$ has proper biharmonic Gauss
map if and only if $n\geq 9$ and
$$
\theta=\frac{1}{2}\arcsin\sqrt{\frac{6n-11\pm\sqrt{4n^2-44n+73}}{2(4n-3)}}.
$$
\end{example}

\subsection*{Isoparametric hypersurface with $\ell = 6$}
In this case, taking into account
\eqref{eq-kalpha}, there exists $\theta\in(0,\pi/6)$ such that
$$
k_1=\cot \theta,\qquad k_2 = \cot\big(\theta +
\frac{\pi}{6}\big)=\frac{\sqrt{3}k_1-1}{k_1+\sqrt{3}},
$$
$$
k_3 = \cot\big(\theta +
\frac{\pi}{3}\big)=\frac{k_1-\sqrt{3}}{1+\sqrt{3}k_1}, \qquad k_4
= \cot\big(\theta + \frac{\pi}{2}\big)=-\frac{1}{k_1},
$$
$$
k_5 = \cot\big(\theta +
\frac{2\pi}{3}\big)=\frac{k_1+\sqrt{3}}{1-\sqrt{3}k_1}, \qquad k_6
= \cot\big(\theta +
\frac{5\pi}{6}\big)=\frac{1+\sqrt{3}k_1}{\sqrt{3}-k_1}.
$$

If $m_i=1$, $i=1,\ldots,6$, then the square of the norm of the
shape operator is
\begin{eqnarray}\label{eq: type 6 A^2}
|\overline{A}|^2&=&\frac{9k_1^{12}+495k_1^8-528k_1^6+495k_1^4+9}{k_1^2(3k_1^4-10k_1^2+3)^2},
\end{eqnarray}
and since $m=6$ the hypercone generated by $\overline{M}$ has
proper biharmonic Gauss map if and only if
$$
|\overline{A}|^2=12.
$$
The last equation, together with \eqref{eq: type 6 A^2}, implies
that $k_1$ is a solution of the equation
$$
x^{12}-12x^{10}+135x^8-216x^6+135x^4-12x^2+1=0,
$$
which has no real roots.

If $m_i=2$, $i=1,\ldots,6$, a similar computation leads to the
conclusion that $k_1$ is a solution of the equation
$$
3x^{12}-45x^{10}+465x^8-766x^6+465x^4-45x^2+3=0,
$$
which has no real roots.

Conclusively, we have
\begin{proposition}\label{cor: isoparam_l=6}
There exist no isoparametric hypersurface of type $6$ whose
associated hypercone has proper biharmonic Gauss map.
\end{proposition}


\begin{thebibliography}{99}

\bibitem{ABR} U.~Abresch. Isoparametric hypersurfaces with four or six
distinct principal curvatures. Necessary conditions on the
multiplicities. {\em Math. Ann.} 264 (1983), 283--302.

\bibitem{AB} A.~Balmu\c s. Biharmonic maps and submanifolds. {\em
Ph.D. Thesis, University of Cagliari, Italy}, January 2008.

\bibitem{ABSMCO1} A.~Balmu\c s, S.~Montaldo, C.~Oniciuc.
Classification results for biharmonic submanifolds in spheres. {\em
Israel J. Math.}, to appear.

\bibitem{RCSMCO1} R.~Caddeo, S.~Montaldo, C.~Oniciuc.
Biharmonic submanifolds of $\s^3$. {\em Internat. J. Math.} 12
(2001), 867--876.

\bibitem{CAR} E.~Cartan. Sur des familles remarquables d'hypersurfaces
isoparam\'{e}triques dans les espaces sphériques. {\em Math. Z.} 45
(1939), 335--367.

\bibitem{EL2} J.~Eells, L.~Lemaire. Selected topics in harmonic maps.
{\em Conf. Board. Math. Sci.} 50 (1983).

\bibitem{ES} J.~Eells, J.H.~Sampson. Harmonic mappings of Riemannian
manifolds. {\em Amer. J. Math.} 86 (1964), 109--160.

\bibitem{OJG} O.J.~Garay. Finite type cones shaped on spherical
submanifolds. {\em Proc. Amer. Math. Soc.} 104 (1988), 868--870.

\bibitem{J1} G.Y.~Jiang. $2$-harmonic isometric immersions between
Riemannian manifolds. {\em Chinese Ann. Math. Ser. A}  7 (1986),
130--144.

\bibitem{Jia} G.Y.~Jiang. 2-harmonic maps and their first and second
variational formulas. {\em Chinese Ann. Math. Ser. A} 7 (1986),
389--402.

\bibitem{LO2} E.~Loubeau,
C.~Oniciuc. The index of biharmonic maps in spheres. {\em
Compositio Math.} 141 (2005), 729--745.

\bibitem{SMCO} S.~Montaldo, C.~Oniciuc. A short survey on biharmonic
maps between Riemannian manifolds. {\em Rev. Un. Mat. Argentina} 47
(2006), 1--22.

\bibitem{MUN} H.F.~M\"{u}nzner. Isoparametrische
Hyperfl\"{a}chen in Sph\"{a}ren. {\em Math. Ann.}  251 (1980),
57--71.

\bibitem{BON} B.~O'Neill. {\em Semi-Riemannian geometry. With applications to
relativity.} Pure and Applied Mathematics, 103. Academic Press,
Inc., New York, 1983.

\bibitem{MO} M.~Obata. The Gauss map of immersions of Riemannian
manifolds in spaces of constant curvature. {\em J.
Differential Geometry} 2 (1968), 217--223.

\bibitem{O2} C.~Oniciuc. Biharmonic maps between Riemannian manifolds. {\em An.
\c Stiin\c t. Univ. Al. I. Cuza Ia\c si. Mat. (N.S.)} 48 (2002),
237--248.

\bibitem{PP} P.~Petersen. {\em Riemannian Geometry}.
Springer-Verlag, 1997.

\bibitem{RV}E.~Ruh, J.~Vilms. The tension field of the
Gauss map. {\em Trans. Amer. Math. Soc.} 149 (1970)
569--573.

\bibitem{JS} J.~Simons. Minimal varieties in Riemannian manifolds.
{\em Ann. Math.} 88 (1968), 62--105.

\bibitem{T1} R.~Takagi. On homogeneous real hypersurfaces in a complex
projective space. {\em Osaka J. Math.} 10 (1973), 495--506.

\end{thebibliography}
\end{document}